\documentclass{article} 
\usepackage{ut-backref} \usepackage{amssymb} \usepackage{amsmath}
\def\covek#1{\mbox{#1}}
\def\diff#1#2{\frac{\partial #1}{\partial #2}}
\def\ldiff#1#2{\partial #1 / \partial #2} 
\def\kaehler{\textup{d}}
\def\KaehlerModule#1#2{\Omega_{#1/#2}} 
\def\rank{\textup{rank}}
\def\grad{\textup{grad}}

\def\card{\textup{\#}} 
\def\height#1{ht\!\left(#1\right)} 
\def\Order{\textup{Order}}
\def\denom{\textup{denom\,}} 
\def\numer{\textup{numer\,}}
\def\lfrac#1#2{#1/#2}
\def\SLP{{\sc slp}} 
\def\System{\Sigma} 

\def\BigO{\mathcal{O}} 
\def\MatMul#1{N(#1)} 
\def\SeriesMul#1{M(#1)}
 
\def\Cindent{\hspace{0.5cm}}

\def\Cwhile{\textbf{while }} 
\def\Cdo{\textbf{do}}

\def\Cif{\textbf{if }}
 
\def\Cthen{\textbf{then }}

\def\Cend{\textbf{end }}

\def\CTest{\textbf{Test}}
\def\CInput{\textbf{Input}\>\>\>} 
\def\COutput{\textbf{Output}\>\>\>}
 
\def\ProcName#1{\small{\sf#1}\normalsize}
\def\Caf{$\gets\ $}
\def\endproof{\hfill\small$\blacksquare$\normalsize\medskip}
\newtheorem{proposition}{Proposition} 
\newtheorem{lemma}{Lemma}
\newtheorem{corollary}{Corollary} 
\newtheorem{definition}{Definition}
 
\newtheorem{theorem}{Theorem}
\title{A probabilistic algorithm to test local algebraic observability
  in polynomial time} \author{Alexandre Sedoglavic \\ Laboratoire
  \textsc{GAGE}, {\'E}cole polytechnique \\ F-91128 Palaiseau, France
  \\ sedoglavic@gage.polytechnique.fr } \date{Preliminary
  version\thanks{This paper is available at~\cite{SedoglavicURL}.  All
  comments are welcome.}~~of \today}
\begin{document}
\maketitle
\begin{abstract}
  The following questions are often encountered in system and control
  theory.  Given an algebraic model of a physical process, which
  variables can be, in theory, deduced from the input-output behavior
  of an experiment?  How many of the remaining variables should we
  assume to be known in order to determine all the others?  These
  questions are parts of the \emph{local algebraic observability}
  problem which is concerned with the existence of a non trivial Lie
  subalgebra of the symmetries of the model letting the inputs and the
  outputs invariant.
 
  We present a \emph{probabilistic seminumerical} algorithm that
  proposes a solution to this problem in \emph{polynomial time}.  A
  bound for the necessary number of arithmetic operations on the
  rational field is presented. This bound is polynomial in the
  \emph{complexity of evaluation} of the model and in the number of
  variables.  Furthermore, we show that the \emph{size} of the
  integers involved in the computations is polynomial in the number of
  variables and in the degree of the differential system.
  
  Last, we estimate the probability of success of our algorithm and we
  present some benchmarks from our Maple implementation.
\end{abstract}
seminumerical algorithm.
\section{Introduction, Notations and Main Result}
\label{sec:introduction}

Local algebraic observability is a structural property of a model and
one of the key-concepts in control theory.  Its earliest definition
goes back to the work of \covek{R.E. Kalman} for the linear case
(see~\cite{Kalman:1961}) and a large literature is devoted to this
subject (see~\cite{HermannKrener:1977, Williamson:1977, Sussmann:1979,
  DiopFliess:1991} and the references therein).  We base our work on
the definition given by \covek{S. Diop}~\& \covek{M. Fliess}
in~\cite{DiopFliess:1991} of the observability for the class of
algebraic systems. \medskip

\begin{figure}[!t]
  \begin{center}
\caption{\label{fig:droso}\small Model for circadian oscillations in
  the Drosophila period protein~\cite{Goldbeter:1995}\normalsize}
\smallskip \newlength{\drosowidth} \setlength{\drosowidth}{\textwidth}
\addtolength{\drosowidth}{-.1\drosowidth}
\fbox{\begin{minipage}{\drosowidth} \smallskip \small$$
    \left\lbrace
      \begin{array}{ccl} \dot{M} & = &
        \frac{v_{s}{K_{I}}^{4}}{{K_{I}}^{4} + {P_{N}}^{4}} -
        \frac{v_{m}M}{K_{m} + M},\\[\smallskipamount]
        \dot{P}_{0} & = & k_{s}M - \frac{V_{1}P_{0}}{K_{1}+P_{0}} +
        \frac{V_{2}P_{1}}{K_{2}+P_{1}},\\[\smallskipamount]
        \dot{P}_{1} & = &
        \frac{V_{1}P_{0}}{K_{1}+P_{0}}+\frac{V_{4}P_{2}}{K_{4} +P_{2}}
        - P_{1}\!\left(\frac{V_{2}}{K_{2}+P_{1}} +
          \frac{V_{3}}{K_{3}+P_{1}}\right)\!,\\[\smallskipamount]
        \dot{P}_{2} & = & \frac{V_{3}P_{1}}{K_{3}+P_{1}}
        -P_{2}\left(\frac{V_{4}}{K_{4}+P_{2}} + k_{1}
          + \frac{v_{d}}{K_{d}+P_{2}}\right) + k_{2}P_{N}, \\[\smallskipamount]
        \dot{P}_{N} & = & k_{1}P_{2} - k_{2}P_{N},\\[\medskipamount]
        y & = & P_{N}.
\end{array} \right.
$$\normalsize\smallskip
\end{minipage}}
\end{center}
\vspace{-\bigskipamount}
\end{figure}
\noindent As in the example of figure~\ref{fig:droso}, such a system is usually
described by means of
\begin{itemize}
\item a vector field, which describes the evolution of \emph{state
    variables} in function of \emph{inputs} and of \emph{parameters};
\item some \emph{outputs} which are algebraic functions of these
  variables.
\end{itemize}

The definition of observability given in~\cite{DiopFliess:1991} relies
on the theory of differential algebra founded by \covek{J.F.
  Ritt}~\cite{Ritt:1966} and is based on the existence of algebraic
relations between the state variables and the successive derivatives
of the inputs and the outputs.

These relations can be considered as an obstruction to the existence
of infinitely many trajectories of the state variables which are
solutions of the vector field and fit the same specified input-output
behavior.  If there are only finitely many such trajectories, the
state variables are said to be locally observable.  \medskip

In order to illustrate this notion, let us consider the \emph{local
  structural identifiability} problem which is a particular case of
the observability problem.  The question is to decide if some unknown
\emph{parameters} of a model are observable considering their
parameters as a special kind of state variables~$\Theta$
satisfying~${\dot{\Theta} = 0}$
(see~\mbox{\cite{RaksanyiLecourtierWalterVenot:1985,
    VajdaGodfreyRabitz:1989, Ollivier:1990:these, GladLjung:1990,
    ChappellGodfreyVajda:1990}}).  If they are not observable, then
infinitely many values of these parameters can fit the same observed
data. Hence, if these parameters have a physical significance, it may
be necessary to change the experimental protocol when possible.  On
the other hand, if the parameters are identifiable, various numerical
approximation methods can be used for their estimation~(see
\cite{Stortelder:1996} and the references therein).  \medskip

We consider the local algebraic observability problem under the
computer algebra standpoint.  The previous studies that enable to test
observability mainly rely on characteristic set or standard bases
computation~\cite{RaksanyiLecourtierWalterVenot:1985,
  Ollivier:1990:these, GladLjung:1990,
  BoulierLazardOllivierPetitot:1995, Hubert:2000} and their complexity
is, at least, exponential in the number of variables and of parameters
(see~\cite{GalloMishra:1991, Sadik:2000:mar}).  Some other techniques,
as the local state-space isomorphism
approach~\cite{VajdaGodfreyRabitz:1989} or the conversion between
characteristic set w.r.t.~different ranking~\cite{Boulier:1999}, can
also be used. The complexities of these methods are not known.

We present a probabilistic polynomial-time algorithm which computes
the set of observable variables of a model and gives the number of non
observable variables which should be assumed to be known in order to
obtain an observable system.  A Maple implementation is available
at~\cite{SedoglavicURL}.  \medskip

\noindent\textbf{Example:} 
Let us consider the use of our algorithm with a model for circadian
oscillations in the Drosophila period protein~\cite{Goldbeter:1995}.
This model is presented in figure~\ref{fig:droso}; there are seventeen
parameters and no input in it.  After one minute of computation, our
Maple implementation gives the following results:
\begin{itemize}
\item the variable~$M$ and the parameters~${\{v_{s}, v_{m}, K_{m},
k_{s}\}}$ are not observable. All the other parameters and variables
are observable;
\item if the non observable variable or only one of the non observable
  parameters are specified, all the variables and parameters of the
  resulting system are observable.
\end{itemize}

Our algorithm certifies that a variable is observable and the answer
for a non observable one is probabilistic with high probability of
success.These results allow us to focus our attention on just four of
the seventeen original parameters. Thus, the search of an
infinitesimal transformation which leaves the output~$y$ and the
vector field invariant is simplified and we find a group of symmetries
generated by~${\{M, v_{s}, v_{m}, K_{m}, k_{s}\} \rightarrow \{\lambda
M, \lambda v_{s}, \lambda v_{m}, \lambda K_{m}, k_{s}/\lambda\}}$.
Hence, there is an infinite number of possible values for non
observable parameters which fit the same specified output~$y$: this
system is certainly unidentifiable.

\subsection{Notations and Main Result}
\label{sec:Notations}
Hereafter, we consider a state-space representation with time
invariant parameters defined by an algebraic system of the following
kind:$$
\stepcounter{equation}
\label{eq:StateSpaceRep} 
\System \quad \left\lbrace
\begin{array}{llll} 
\dot{\Theta} & = & 0, \\ 
\dot{X} & = & F(X,\Theta,U), & \hspace{1cm} (\theequation.1) 
\\[\smallskipamount] 
Y & = & G(X,\Theta,U). & \hspace{1cm} (\theequation.2) 
\end{array} \right.
$$
Big letters stand for vector-valued objects and we suppose that
there are:
\begin{itemize}
\item $\ell$ parameters~${\Theta:=(\theta_{1},\dots,\theta_{\ell})}$
\item $n$ state variables~${X:=(x_{1},\dots,x_{n})}$;
\item $r$ input variables~${U:=(u_{1},\dots,u_{r})}$;
\item $m$ outputs variables~${Y:=(y_{1},\dots,y_{m})}$ with~${m \leq
    n}$.
\end{itemize} 
The letter $\dot{X}$ stands for the derivatives of the state
variables~${(\dot{x}_{1}, \dots, \dot{x}_{n})}$ and~$F$ (resp.~$G$)
represents~$n$ (resp.~$m$) rational fractions
in~$\mathbb{Q}(X,\Theta,U)$ which are denoted by~$(f_{1},\dots,f_{n})$
\big(resp.  $(g_{1}, \dots,g_{m})$\big).  The letter~$d$ (resp.~$h$)
represents a bound on the degree (resp. size of the coefficients) of
the numerators and denominators of the~$f_{i}$'s and~$g_{i}$'s.

Hereafter, we use a common encoding where the expression~${e := x^5}$
is represented as a sequence of instructions:~${t_{1}:= x, t_{2}:=
  {t_{1}}^{2}, t_{3}:={t_{2}}^{2}, e:=t_{3} t_{1}}$.

Hence, the system~$\System$ is represented by a \emph{straight-line
  program} without division which computes its numerators and
denominators and requires~$L$ arithmetic operations (see
Section~\ref{sec:Complexity} and~{\S}~4
in~\cite{BurgisserClausenShokrollahi:1997}). \smallskip

The following theorem is the main result of this paper.
\begin{theorem}
  \label{th:ComplexityOfLAOT}

  Let~$\System$ be a differential system as described in
  Section~\ref{sec:Notations}.  There exists a probabilistic algorithm
  which determines the set of observable variables of~$\System$ and
  gives the number of non observable variables which should be assumed
  to be known in order to obtain an observable system.
  
  The arithmetic complexity of this algorithm is bounded
  by$$\BigO\Biggl(\! M(\nu)\Bigl(\MatMul{n+\ell} + (n+m)L\Bigr) +
  (n+\ell+1)\MatMul{n+\ell}\frac{m\nu}{n+\ell} \!\Biggr)$$
  with~$\SeriesMul{\nu}$ \big(resp.  $\MatMul{\nu}$\big) the cost of power
  series multiplication at order~${\nu+1}$ (resp.~${\nu \times \nu}$ matrix
  multiplication) where~${\nu \leq n+\ell}$.
  
  Let~$\mu $ be an arbitrary positive integer,~$D$
  be~${4(n+\ell)^{2}(n+m)d}$ and \small$$
  D^{\prime} := \big(2
  \ln(n+\ell+r+1) + \ln \mu D\big)D + 4(n+\ell)^{2}\big( (n+m)h + \ln 2nD
  \big).
  $$\normalsize If the computations are done modulo a prime number~${p
    > 2D^{\prime}\mu }$ then the probability of a correct answer is at
  least~${(1-1/\mu )}^{2}$.  \end{theorem}
 
For the model presented in figure~\ref{fig:droso}, the significant
terms of our complexity statement are~${L=91, n=5, \ell=17, m=1, d=6,
  h=1, \nu=n+\ell}$. The choice of~${\mu =3000}$ leads to a probability
of success around~$.9993$ and the computations are done
modulo~$10859887151$.  These computations take~$10$ seconds on a PC
Pentium~III~(633 Mhz) provided by the UMS MEDICIS~\cite{MEDICISURL}.

\paragraph{Outline of the paper:}
In the next section, we recall some basic definitions of differential
algebra and the definition of algebraic observability used by
\covek{S. Diop}~\& \covek{M. Fliess}
in~\cite{DiopFliess:1991}. Furthermore, we describe the relationship
between this framework and the approach of \covek{H. Pohjanpalo}
in~\cite{Pohjanpalo:1978}.  Then, we present an algebraic jacobian
matrix which is derived from the theory of K{\"a}hler differentials
and used in the local algebraic observability test. \smallskip

In the second part of this paper, we present some new results.  In
Section~\ref{sec:Algorithm}, we show how to compute some
specializations of this matrix using power series expansion of the
output and we estimate the related arithmetic complexity.  Then, we
study the behavior of the integers involved in the computations and we
precise the probabilistic aspect.  In conclusion, we present some
benchmarks.

\section{Differential Algebra and Observability}
\label{sec:setting}
  
Differential algebra, founded by \covek{J.F. Ritt}, is an appropriate
framework for the definition of algebraic observability introduced by
\covek{S. Diop}~\& \covek{M. Fliess} in~\cite{DiopFliess:1991}.  For
more details on differential algebra, we refer to~\cite{Ritt:1966}
and~\cite{Kolchin:1973}; nevertheless, we recall briefly some
necessary notions.

\subsection{Differential Algebraic setting}
\label{sec:DiffAlg}

Let us denote by~$k$ a base field of characteristic zero.  The
differential algebra~${k\{U\}}$ is the~$k$-algebra of multivariate
polynomials defined by the infinite set of indeterminates~${\{U^{(j)}
  |\, \forall j \in \mathbb{N}^{\star}\}}$ and equipped with a
derivation~$\delta$ such that~${\delta u^{(i)} = u^{(i+1)}}$. Its
differential fraction field is denoted by~${k \langle U \rangle}$.
\medskip

\noindent\textbf{Hypotheses:}
\label{sec:Hypotheses} 
The inputs~$U$ and all their derivatives are assumed to be
independent.  Furthermore, we consider non singular solutions
of~$\System$; thus, we assume that we work in an open set where the
denominators present in~$\System$ do not vanish. These hypotheses
represent practically all the encountered systems. \medskip

\subsection{Local Algebraic Observability}
\label{sec:ObservabilityDefinition}

Following the interpretation due to \covek{M. Fliess} of some
algebraic control theory problems~\cite{Fliess:1989}, we consider the
differential field~${\mathcal{K}:=k\langle U \rangle (X,\Theta)}$
equipped with the following formal Lie derivation:$$
\mathcal{L} :=
\diff{}{t} + \sum^{n}_{i=1} f_{i}\diff{\;\;}{x_{i}} + \sum_{j \in
  \mathbb{N}} \sum_{u\in U} u^{(j+1)} \diff{\quad\;}{u^{(j)}}.
$$
This derivation is associated with the vector field defined by the
equations~(1.1). Hereafter, we denote~${(\mathcal{L}f_{1}, \dots,
\mathcal{L}f_{n})}$ by~$\mathcal{L}F$ and~${\underbrace{\mathcal{L}
\circ \dots \circ \mathcal{L}}_{j \textup{ times}}}$
by~$\mathcal{L}^{j}$.  

\noindent Hence, the outputs~$G(X,\Theta,U)$ are denoted
by~$Y$ and~${Y^{(j)} =\mathcal{L}^{j}G(X,\Theta,U)}$.
\begin{definition}[\cite{GladLjung:1990, DiopFliess:1991}]
\label{def:FieldExtensionFormulation}
An element~$z$ in~$\mathcal{K}$ is locally algebraically observable
with respect to inputs and outputs if it is algebraic over~$k\langle
U,Y \rangle$.  Thus, the system~$\System$ is locally observable if the
field extension~${k\langle U,Y \rangle \hookrightarrow \mathcal{K}}$
is purely algebraic.
\end{definition}

Let us illustrate this definition with the following example:$$
\left\lbrace
\begin{array}{lcl}
  \dot{x}_{3} & = & \theta x_{1},\\ 
  \dot{x}_{2} & = & \lfrac{x_{3}}{x_{2}},\\ 
  \dot{x}_{1} & = & \lfrac{x_{2}}{x_{1}},\\[\smallskipamount]
  y & = & x_{1}.
\end{array}
\right.
$$
By successive differentiations of the output, we obtain the
following differential relations:$$
\begin{array}{c}
y -x_{1}, \quad y\dot{y} - x_{2}, \quad y\dot{y}(\dot{y}^{2} + y\ddot{y}) -
x_{3} \\[\smallskipamount] {(\dot{y}^{2} + y\ddot{y})^{2} +
y\dot{y}\big(3\dot{y}\ddot{y} + yy^{(3)}\big) - \theta y}.
\end{array}
$$
Thus, the parameter and the variables are observable according to
Definition~\ref{def:FieldExtensionFormulation}.  Furthermore, as these
algebraic relations define a unique solution, the parameter and the
variables are said to be \emph{globally} algebraically
observable~\cite{GladLjung:1990, Ollivier:1990:these,
  ChappellGodfreyVajda:1990}.

These relations depend generically of high order derivatives of the
output and thus, they are not of a great practical interest for
parameter estimation.  As we focus our attention on local
observability, we are going to avoid their computation. \medskip

\noindent\textbf{Convention:}\label{sec:Convention} We wish to test
observability for the parameters~$\Theta$ and/or state variables~$X$.
Thus, we present the algorithm in the most general case (observability
of parameters and state variables) and we do not describe the
restriction to one case or the other. \medskip

Definition~\ref{def:FieldExtensionFormulation} implies that local
algebraic observability is related to the transcendence degree of the
field extension~${k\langle U, Y \rangle \hookrightarrow \mathcal{K}}$.
So, this property can be tested by a rank computation using K{\"a}hler
differentials (see~Section~\ref{sec:RankCondition}).  As noticed
in~\cite{DiopFliess:1991}, this approach leads to a condition which is
the formal counterpart of the \covek{R. Hermann}~\& \covek{A. Krener}
rank condition in the differential geometric point of
view~\cite{HermannKrener:1977}. \smallskip

Furthermore, the transcendence degree of the field
extension~${k\langle U, Y \rangle \hookrightarrow \mathcal{K}}$ is the
number of non observable variables which should be assumed to be known
in order to obtain an observable system. Thus,
Theorem~\ref{th:ComplexityOfLAOT} is based on the study of this field
extension.

\subsection{A Description of~${k\langle U, Y \rangle \hookrightarrow \mathcal{K}}$}
\label{sec:Generators}
Let us denote by~$\Phi(X,\Theta,U,t)$ the formal power series with
coefficients in~$\mathcal{K}$ such that~${\Phi(X,\Theta,U,0):=X}$
and~${\dot{\Phi} = F(\Phi,\Theta,U)}$, we have:$$
\Phi(X,\Theta,U,t)=
X + \sum_{j \in \mathbb{N}^{\star}} \mathcal{L}^{j} F(X,\Theta,U) \;
\frac{t^{j}}{j!}.$$
Furthermore, let us define the formal power
series~$Y(X,\Theta,U,t)$ with coefficients in~$\mathcal{K}$ such
that~${Y(X,\Theta,U,t) := G\big(
  \Phi(X,\Theta,U,t),\Theta,U,t\big)}$:~\begin{equation}
  \label{eq:OutputSeries} Y(X,\Theta,U,t)  = G(X,\Theta,U) + \sum_{j
    \in \mathbb{N}^{\star}} \mathcal{L}^{j} G(X,\Theta,U) \;
  \frac{t^{j}}{j!}.
\end{equation} 
We recall that these expressions are
vector-valued~${\big(Y=(y_{1},\dots,y_{m})\big)}$. \medskip

In~\cite{Pohjanpalo:1978}, \covek{H. Pohjanpalo} already considers the
coefficients of the power series~$Y(X,\Theta,U,t)$ in order to test
identifiability.  In~\cite{DiopFliess:1991}, the authors prove that a
finite number of these coefficients are necessary to {\it describe}
the field extension~${k\langle U, Y \rangle \hookrightarrow
\mathcal{K}}$. But in these two papers the necessary order of
derivation is not bounded. \smallskip

This can be done using the differential algebra point of view
(see~{\S}~4 in~\cite{Sadik:2000:mar} for a general statement).  The
following proposition summarizes these results in a field extension
framework.
\begin{proposition}
\label{prop:Generators}
The field~$k \langle U, Y \rangle$ is isomorphic to~$k\langle U
\rangle \big(Y,\dots,Y^{(n+\ell+1)}\big)$ and algebraic over~$k\langle U
\rangle \big(Y,\dots,Y^{(n+\ell)}\big)$.
\end{proposition}
\textbf{Proof:} The transcendence degree of~${k\langle U \rangle
  \hookrightarrow \mathcal{K}}$ is equal to~${n+\ell}$.  Hence, the
transcendence degree of~${k\langle U \rangle \hookrightarrow k \langle
  U, Y \rangle}$ is bounded by~${n+\ell}$.  It means that,
for~${i=1,\dots,m}$, there is an algebraic
relation~${q_{i}\big(y_{i},\dots,{y_{i}}^{(n+\ell)}\big)=0}$ and the
derivative~${y_{i}}^{(n+\ell+1)}$ is a rational function
of~$y_{i},\dots,{y_{i}}^{(n+\ell)}$ with coefficients in~$k \langle U
\rangle$.  This proves that~$k \langle U, Y \rangle$ is isomorphic
to~$k\langle U \rangle \big(Y,\dots,Y^{(n+\ell)}\big)$. \endproof

If there is more than a single output, the necessary order of
derivation can be smaller than~${n+\ell}$ and it is denoted by~$\nu$.
This index of differentiation is a natural measure of the complexity
of our algorithm (see Section~\ref{sec:Complexity}) and
generically~${\nu=(n+\ell)/m}$. Hereafter, we take~$\nu$ equal
to~${n+\ell}$ as in Theorem~\ref{th:ComplexityOfLAOT}. \smallskip

In the above proof, following the hypotheses of
Section~\ref{sec:Hypotheses}, we assumed that the independent input
variables~$U$ and all their derivatives were in the base field.
Furthermore, we showed that we just need the first~${n+\ell}$
derivatives of the output equations.  In order to simplify the
presentation in the next section, we assume that the base field
is~${\bar{k}:=k\big(U,Y,\dots,U^{(n+\ell)},Y^{(n+\ell)}\big)}$.

We present now the properties of the module of K{\"a}hler
differentials which are used to compute the transcendence degree
of~${\bar{k} \hookrightarrow \bar{k}(X, \Theta)}$ in practice.

\subsection{Rank Conditions}
\label{sec:RankCondition}
If~$S \hookrightarrow T$ is a field extension, we use the
notation~$\KaehlerModule{T}{S}$ for the $T$-vector space which is the
cokernel of the jacobian matrix~$\ldiff{(\mathcal{L}^{i} G )_{0 \leq i
\leq \nu}} {(X,\Theta)}$ and~${\kaehler z}$ stands for the image
of~${z \in T}$ in this vector space (see {\S}~16
in~\cite{Eisenbud:1994} for standard definition
and~\cite{Johnson:1969} for construction in differential algebra).  We
recall the following result:
\begin{theorem}[{\S}~16 in~\cite{Eisenbud:1994}]
  \label{thm:KaelherModuleGenerators} Let us consider~$S$ a field of characteristic zero and~$T$ a finitely generated
  field extension of~$S$. If~${\{ x_{\lambda}\}_{\lambda \in \Lambda}
    \subset T}$ is a collection of elements, then~${\{\kaehler
    x_{\lambda}\}_{\lambda \in \Lambda}}$ is a basis
  of~$\KaehlerModule{T}{S}$ as a vector space over~$T$ iff the~$\{
  x_{\lambda}\}_{\lambda \in \Lambda}$ form a transcendence basis
  of~$T$ over~$S$.
\end{theorem}
Our algorithm is based on the following straightforward consequences
of this theorem.
\begin{corollary} 
\label{cor:RankCondition} 
If~$\phi$ is the transcendence degree of the field extension~${\bar{k}
  \hookrightarrow \bar{k}(X, \Theta)}$ then we have the equality$$\phi
  = (n+\ell) - \rank_{\bar{k}(X, \Theta)}\!\left(
  \ldiff{\!\left(\mathcal{L}^{i} G \right)_{0 \leq i \leq \nu}}
  {(X,\Theta)} \right)\!.$$
\noindent Furthermore, If the rank of the jacobian
submatrix~${\ldiff{(\mathcal{L}^{j}G)_{0 \leq j \leq \nu}}
{(X\!\setminus\!\{x_{i}\},\Theta)}}$
(resp.~${\ldiff{(\mathcal{L}^{j}G)_{0 \leq j \leq \nu}}
{(X,\Theta\!\setminus\!\{\theta_{i}\})}}$) is equal
to~${n+\ell-\phi}$, then the transcendence degree of the field
extension~$\bar{k}\hookrightarrow \bar{k}(x_{i})$
(resp.~$\bar{k}\hookrightarrow \bar{k}(\theta_{i})$) is equal to zero
and the variable~$x_{i}$ (resp. the parameter~$\theta_{i}$) is
observable.
\end{corollary} 
The computation of~$\phi$ is mainly based on the construction and the
evaluations of a \emph{straight-line program} which allows to compute
the power series expansion of~$Y(X,\Theta,U,t)$.  We present the
necessary notions in the next section.

\subsection{Data Encoding and Complexity Model}
\label{sec:SLP}
The above results can be expressed considering a polynomial~$f$ as an
element of a vector space; hereafter, we consider an algebraic
expression as a function. \medskip

This classical point of view in numerical analysis is also used in
computer algebra for complexity statements or practical algorithms
(\mbox{see~\cite{ GiustiLecerfSalvy:1999, GathenGerhard:1999,
Schost:2000, Sedoglavic:2000:jan}} and the references therein).  We
refer to Chapter~4 of~\cite{BurgisserClausenShokrollahi:1997} for more
details about this model of computation.
\begin{definition}
  Let~${A := \{a_{1},\dots,a_{j}\}}$ be a finite set of variables. A
  straight-line program is a sequence of assignments~${b_{i}
    \leftarrow b^{\prime} \circ_{i} b^{\prime\prime}}$ where
  ${\circ_{i} \in \{+,-,\times ,{\div}\}}$ and where~${\{b^{\prime}
    ,b^{\prime\prime}\} \subset \bigcup^{i-1}_{j=1} \{ b_{j} \} \cup A
    \cup k}$.  Its complexity of evaluation is measured by its
  length~$L$, which is the number of its arithmetic operations.
  Hereafter, we use the abbreviation \SLP\ for straight-line program.
\end{definition}

As a \SLP\ representing a rational expression~${f \in k(a_{1}, \dots,
  a_{j})}$ is a program which computes the value of~$f$ from any
  values of the base field such that every division of the program is
  possible.  Furthermore, it is possible to determine a \SLP\
  representing the gradient of~$f$.  The following constructive
  results allows us to handle these two aspects.
\begin{theorem}[\covek{W. Baur}~\& \covek{V. Strassen}~\cite{BaurStrassen:1983}]
  \label{prop:SLPProp} Let us consider a \SLP\ computing the value of
  a rational expression~$f$ in a point of the base field and let us
  denote by~$L_{f}$ its complexity of evaluation.
 
  One can construct a \SLP\ of length~$5L_{f}$ which computes the
  value of~$\grad (f)$.
\end{theorem}
Furthermore, one can construct a \SLP\ of length~$4L_{f}$ which
computes two polynomials~$f_{1}$ and~$f_{2}$ such that~${f =
f_{1}/f_{2}}$. \medskip

Following our presentation, one can construct formally all the
expressions introduced in Sections~\ref{sec:Generators}
and~\ref{sec:RankCondition} with its favourite computer algebra
system.  

But, let us recall that, in order to compute the formal
expressions~$\mathcal{L}^{\nu}G$ and the associated jacobian matrix,
one has to differentiate~$\nu$ times the output equations~(1.2).  As
explained in~\cite{Kaltofen:1993}, the arithmetic complexity of
computing multiple partial derivatives is likely exponential in~$\nu$.
If the evaluation complexity of the output equations~(1.2) is~$L$, by
Theorem~\ref{prop:SLPProp}, the computation of~$\mathcal{L}^{\nu} G$
requires at least~$(5m)^{\nu}L$ arithmetic operations.  

Thus, this strategy cannot lead to a polynomial time algorithm.
\medskip

The rank computations defined in the previous section are also
cumbersome because they are mainly performed on the
field~$\bar{k}(X,\Theta)$.  Nevertheless, in order to determine~$\phi$
efficiently, the variables~$X$,~$\Theta$ and~$U$ can be specialized to
some generic values in the jacobian matrix and so, its generic rank
can be computed numerically with high probability of success (see
Section~\ref{sec:ProbabilisticAspect}). \medskip

Thus, the main problem is to avoid the formal computation
of~$(\mathcal{L}^{i} G)_{0 \leq i \leq \nu}$.  In fact, our strategy
is to specialize a linearized system derived form~$\System$ first and
to recover the value of~$\phi$ just using numerical computations on a
finite field.

\section{A Probabilistic Polynomial-Time Algorithm}
\label{sec:Algorithm}

In Section~\ref{sec:JacobianMatrixComputation}, we present the
\emph{linear variational system} derived from~$\System$ which allows
us to compute directly the jacobian
matrix~$\ldiff{(\mathcal{L}^{i}G)_{0 \leq i \leq \nu}} {(X, \Theta)}$
with~$X,\Theta$ and~$U$ specialized on some given values.  

Then, we show how this matrix can be determined in polynomial time and
we give an estimation of the arithmetic complexity of our
algorithm. \medskip

The purpose of the Sections~\ref{subsubsec:Bound}
and~\ref{sec:ProbabilisticAspect} is to study the growth of the
integers involved in the computations and to estimate the probability
of success of our algorithm.

\subsection{Variational System Derived From~$\System$}
\label{sec:JacobianMatrixComputation}
As shown in Section~\ref{sec:RankCondition}, our goal is to compute
the generic rank of the jacobian matrix~$\ldiff{(\mathcal{L}^{i} G)_{0
    \leq i \leq \nu}} {(X,\Theta)}$.  Using
relation~(\ref{eq:OutputSeries}), we conclude that:$$
\diff{(\mathcal{L}^{j}G)_{0 \leq j \leq \nu}}{(X,\Theta)\;\quad} =
\diff{\big(\textup{coeffs}\big(Y(t)\big)\big)}{(X,\Theta)} =
\textup{\ProcName{coeffs}}\left( \diff{G}{X} \diff{\Phi}{X\;} ,
  \diff{G}{X} \diff{\Phi}{\Theta} + \diff{G}{\Theta}\!\right)\!.$$
The
above equalities leads to the following relation:
\begin{equation}
   \label{eq:Output} \diff{(\mathcal{L}^{j}G)_{0 \leq j \leq
\nu}}{(X,\Theta)\;\quad} = \textup{\ProcName{coeffs}}
\left(\nabla Y\left(\Phi, \diff{\Phi}{X},
\diff{\Phi}{\Theta}\right)\!, t^{j}, j=0,\dots,\nu\right)\!,
\end{equation}
where~$\nabla Y$ denote the following~${n\times (n+\ell)}$ matrix
represented by a \SLP:$$
\nabla Y
\big(\Phi,\Gamma,\Lambda,\Theta,U\big) := \left( \diff{G}{X} \Gamma,
  \diff{G}{X} \Lambda + \diff{G}{\Theta}\!\right)
\!\big(\Phi,\Gamma,\Lambda,\Theta,U\big)\!.
$$
Hence, we have to determine the first~${\nu=n+\ell}$ terms of the power
series expansion of~$\Phi(X,\Theta,U,t)$,~${\Gamma(X,\Theta,U,t)
  :=\ldiff{\Phi}{X}}$ and~${\Lambda(X,\Theta,U,t) :=
  \ldiff{\Phi}{\Theta}}$. \medskip

Let us denote by~${P(\dot{X},X,\Theta,U) = 0}$, the numerators of the
rational relations~${\dot{X} - F(X,\Theta,U) = 0}$ and let us consider
the following expressions:
\begin{equation}
\label{eq:System}
\nabla P 
\left\lbrace
   \begin{array}{cl}
      P(\dot{X},X,\Theta,U), & (\ref{eq:System}.1) 
      \\[\smallskipamount]
      \diff{P}{\dot{X}}\big(X,\Theta,U\big) \dot{\Gamma} +
      \diff{P}{X}\big(\dot{X},X,\Theta,U\big) \Gamma,&(\ref{eq:System}.2)
      \\[\smallskipamount]
      \diff{P}{\dot{X}}\big(X,\Theta,U\big) \dot{\Lambda} +
      \diff{P}{X}\big(\dot{X},X,\Theta,U\big) \Lambda +
      \diff{P}{\Theta}\big(\dot{X},X,\Theta,U\big). & (\ref{eq:System}.3)
\end{array} \right.
\end{equation}
The power series $\Phi(X,\Theta,U,t)$,~$\Gamma(X,\Theta,U,t)$
and~$\Lambda(X,\Theta,U,t)$ are solutions of the system of ordinary
differential equations~${\nabla P(\Phi,\Gamma,\Lambda,\Theta,U)=0}$
with initial conditions~${\Gamma(X,\Theta,U,0):=\mathrm{Id}_{n \times n}}$
and~${\Lambda(X,\Theta,U,0):=0_{n \times \ell}}$. \medskip

\noindent\textbf{Commentary:} We have already noticed that one can
compute symbolically the expression of the formal jacobian
matrix~$\ldiff{(\mathcal{L}^{i}G)_{0 \leq i \leq \nu}} {(X, \Theta)}$.
The rank computations described in Corollary~\ref{cor:RankCondition}
are sufficient to conclude.  

Furthermore, if~$X,\Theta$ and~$U$ are specialized on some random
values, these computations can be performed numerically with high
probability of success.  We summarize this possible strategy in the
upper horizontal and the right vertical arrow of the following
diagram:
\begin{center}
  \def\PictureWidth{300} \def\PictureHeight{90}
  \begin{picture}(\PictureWidth,\PictureHeight) 
     \newsavebox{\specmorp}
     \savebox{\specmorp}(0,0){\scriptsize
        $\left\lbrace 
        \begin{array}{ccl} 
           X&\rightarrow&X_{0} \in \mathbb{Z}^{n}, \\
           \Theta&\rightarrow&\widetilde{\Theta} \in \mathbb{Z}^{\ell}, \\
           U&\rightarrow&\widetilde{U} \in (\mathbb{Z}[t])^{r}.  
        \end{array} \right.$\normalsize}
     \newsavebox{\jacmat} 
     \savebox{\jacmat}(0,0){
        $\left(\diff{\left(\mathcal{L}^{i}G\right)_{0\leq i\leq\nu}}
        {(X,\Theta)\quad}\right)$} 
     \thinlines 
     \put(125,83){\makebox(0,0){formal computation}}
     \put(125,80){\makebox(0,0){\vector(1,0){125}}} 
     \put(230,75){\usebox{\jacmat}} 
     \put(25,75){\makebox(0,0){$\System$}} 
     \put(225,50){\makebox(0,0){\vector(0,-1){30}}}
     \put(280,47){\usebox{\specmorp}}
     \thicklines 
     \put(25,50){\makebox(0,0){\vector(0,-1){30}}}
     \put(25,15){\makebox(0,0)
        {$\nabla P$}}
     \put(125,25){\makebox(0,0){numerical computation on~$\mathbb{Q}$}}
     \put(125,20){\makebox(0,0){\vector(1,0){125}}}
     \put(230,15){\usebox{\jacmat}} 
     \put(285,15){ \makebox(0,0)
        {$\big(X_{0},\widetilde{\Theta},\widetilde{U}\big)$}}
\end{picture}
\end{center}
As the symbolic computation of the jacobian matrix is cumbersome, we
specialize the parameters on some random integers~$\widetilde{\Theta}$
and the inputs~$U$ on the power series~$\widetilde{U}$ which are
truncated at order~${n+\ell+1}$ with random integer coefficients.
Then, we solve the associated system~$\nabla P$ for some integer
initial conditions~$X_{0}$ and we compute the
specialization~$\ldiff{(\mathcal{L}^{i}G)_{0 \leq i \leq \nu}} {(X,
  \Theta)}(X_{0},\widetilde{\Theta})$ with~$\nabla Y$.  This approach
is summarized by the left vertical and the lower horizontal arrow.  We
present an algorithm which relies on this standpoint and we give in
Section~\ref{sec:ProbabilisticAspect} its probability of success.
\medskip

The hypothesis~${\ldiff{P}{\dot{X}} \neq 0}$ assumed in
Section~\ref{sec:Hypotheses} ensures that the differential system
${\nabla P\big(\Phi, \Gamma,
  \Lambda,\widetilde{\Theta},\widetilde{U}\big) = 0}$ admits an unique
formal solution~\cite{DenefLipshitz:1984} which can be computed with
the following Newton operator.

\subsection{A Quadratic Newton Operator}
\label{sec:Newton}
The aim of this section is to present the Newton operator used in our
algorithm.  In~\cite{Geddes:1979, BrentKung:1978}, the authors show
that such an operator is quadratic. We sketch its construction and
neglect the technical details for the sake of simplicity. \medskip

We recall that we work with vector-valued expressions. Thus, the
expression~(\ref{eq:System}.1)
\big(resp.~(\ref{eq:System}.2),~(\ref{eq:System}.3)\big) represents
a~${n \times 1}$ (resp.~${n\times n,\ n\times \ell}$) matrix.

The Theorem~\ref{prop:SLPProp} allows to construct, from a \SLP\ of
length~$L$ which encodes~$\System$, another \SLP\ of
length~$\BigO\big(N(n+\ell)+nL\big)$ which encodes the system~$\nabla
P$. For some given series~$\Phi,\Gamma$ and~$\Lambda$, this \SLP\ 
computes the following~${n\times (1+n+\ell)}$ matrix:$$
\left(
\begin{array}{ccl}
  p_{1}(\dot{\Phi},\Phi,\widetilde{\Theta},\widetilde{U}) &

\diff{P}{\dot{X}}\big(\Phi,\widetilde{\Theta},\widetilde{U}\big) \dot{\Gamma} &
\diff{P}{\dot{X}}\big(\Phi,\widetilde{\Theta},\widetilde{U}\big) \dot{\Lambda}
\hfill + \\ 

\vdots & + &\diff{P}{X}\big(\dot{\Phi}, \Phi, \widetilde{\Theta}, 
\widetilde{U}\big) \Lambda \; + \\[\medskipamount] 

p_{n}(\dot{\Phi},\Phi,\widetilde{\Theta},\widetilde{U}) &
\diff{P}{X}\big(\dot{\Phi},\Phi,\widetilde{\Theta},\widetilde{U}\big) \Gamma &
\diff{P}{\Theta}\big(\dot{\Phi},\Phi,\widetilde{\Theta},\widetilde{U}\big) \\

\end{array}
\right)\!.
$$

Let us represent~$\Phi(X,\Theta,U,t)$
\big(resp.~$\Lambda(X,\Theta,U,t),\ \Gamma(X,\Theta,U,t)$\big)
mod~$t^{2^{j}}$ by~$\Phi_{j}$ \big(resp.~$\Lambda_{j},\ 
\Gamma_{j}$\big) and denote the correction term:$${\big(
  \Phi(X,\Theta,U,t) - \Phi_{j}, \Gamma(X,\Theta,U,t) - \Gamma_{j},
  \Lambda(X,\Theta,U,t) - \Lambda_{j} \big)} \bmod {t^{2}}^{j+1}
\textup{ by } E_{j+1}.$$
\smallskip

As usually, we construct our Newton operator from the Taylor series
expansion of the function~$\nabla P$. This yields the following
relations: \small$$
\nabla P
\big(\Phi,\Gamma,\Lambda\big)(X,\Theta,U,t) =\!  \nabla P\big(
\Phi_{j}, \Gamma_{j}, \Lambda_{j}\big) + \diff{\;\nabla P\quad}
{(\dot{X},\dot{\Gamma},\dot{\Lambda})}{} \dot{E}_{j+1} +
\diff{\;\nabla P\quad} {(X,\Gamma,\Lambda)} {}E_{j+1} + \dots\! = 0.
$$\normalsize The remaining terms are of order in~$t$ greater
than~$2^{j+1}$. Thus, they are not necessary for the computation
of~$E_{j}$. \medskip

\noindent\textbf{Computational strategy:} 
we consider~$\Phi$ as a variable in the first column of~$\nabla P$ and
as a constant in the others. Thus, we have the following
relations:\small
$$
\diff{\;\nabla P\quad}{(\dot{X}, \dot{\Gamma}, \dot{\Lambda})}{} =
\left(\diff{P}{\dot{X}}{} , \diff{P}{\dot{X}}{},\diff{P}{\dot{X}}{}
\right)\!, \quad \diff{\;\nabla P\quad}{(X,\Gamma,\Lambda)}{} =
\left(\diff{P}{X}{},\diff{P}{X}{},\diff{P}{X}{}\right)\!.
$$\normalsize
\noindent\textbf{Consequence of our computational strategy:} The
above hypothesis induces a \emph{shift} between the order of correct
coefficients of~$\Lambda_{j}$,~$\Gamma_{j}$ and~$\Phi_{j}$.  In
fact,~$\Lambda_{j}$ and~$\Gamma_{j}$ are correct
modulo~${t^{2}}^{j-1}$.  Thus, we need to stop the following operator
with~${{j+1} = \ln_{2} (n+\ell+1)}$ and to repeat one more time the
last resolution at the same order. \bigskip

\noindent\textbf{Newton operator:}
The above hypothesis leads to a Newton operator based on the
resolution of the following system of linear ordinary differential
equations: \small
\begin{equation}
  \label{eq:NewtonOperator}
\!\!\!\diff{P}{\dot{X}}\big(\Phi_{j},\widetilde{\Theta},\widetilde{U}\big)
\dot{E}_{j+1}\! +\!
\diff{P}{X}\big(\dot{\Phi}_{j},\Phi_{j},\widetilde{\Theta},\widetilde{U}\big)
E_{j+1}\! +\! \nabla P \big(\Phi_{j}, \Gamma_{j},\Lambda_{j}
,\widetilde{\Theta},\widetilde{U} \big)\!= 0
\bmod t^{2^{j+1}}\!
\end{equation}
\normalsize From the initial conditions~${\Phi_{0}\in
  \mathbb{Z}^{n}}$,~${\Gamma_{0}:=\mathrm{Id}_{n\times n}}$
and~${\Lambda_{0}:=0_{n \times \ell}}$, this system is solved iteratively
for~${{j+1} = 1,\dots,\ln_{2} (n+\ell+1)}$ using the recurrence
relations~${\big(\Phi_{j+1}, \Gamma_{j+1}, \Lambda_{j+1}\big) =
  \big(\Phi_{j}, \Gamma_{j}, \Lambda_{j}\big) + E_{j+1}}$. \medskip

The resolution of the linear ordinary differential
system~(\ref{eq:NewtonOperator}) relies on the method of integrating
factors.  First, we consider the Homogeneous system \small$$
\diff{P}{\dot{X}} \big(\Phi_{j},\widetilde{\Theta},\widetilde{U}\big)
\dot{W}_{j} + \diff{P}{X}
\big(\dot{\Phi}_{j},\Phi_{j},\widetilde{\Theta},\widetilde{U}\big)
W_{j} = 0 \bmod t^{2^{j+1}}
$$\normalsize where~$W_{j}$ denote a~${n \times n}$ unknown matrix which
coefficients are series truncated at order~$2^{j}$.  The main trick is
common in power series manipulation, we consider matrices with
coefficients in a series ring as series with coefficients in a matrix
ring. For example, we have~${A \bmod t^{2^{j+1}} = A_{0} + A_{1}t +
  \dots + A_{2^{j}} t^{2^{j}}}$ where the~$A_{i}$'s are matrices with
coefficients in the rational field. \medskip

Thus, the product, the exponential and, if~$A_{0}$ is invertible, the
inverse of matrices with coefficients in a series ring can be computed
at precision~$j$ with the classical Newton operator (see~4.7
in~\cite{Knuth:1998} and {\S}~5.2 in~\cite{BrentKung:1978} for more
details).  For example, if~$A_{0}$ is invertible and~$B_{j}$ denotes
the inverse of~$A$ at order~$t^{2^{j}}$, we have~$B_{j+1} = 2B_{j} -
B_{j}AB_{j}$.\medskip

Furthermore, it is a basic fact from the theory of linear ordinary
system that if~${A\dot{W}+A^{\prime}W=0}$ and~$A$ is invertible
then~${W= \exp\!\big(\!\int A^{-1}A^{\prime}\big)}$ is a matricial
solution of this system.  Hence, the above homogeneous system can be
solved at precision~$j$ by a procedure called
\ProcName{HomogeneousResolution} in figure~\ref{fig:Algorithm}.
\medskip

With the same tools, one can check that the following formal
expression deduced from the formula for variation of
constants~\small$$
W^{-1} \int \left(W
  \left(\diff{P}{\dot{X}}{}\right)^{-1}\!  \nabla P\right) \big(
\Phi_{j}, \Gamma_{j},
\Lambda_{j},\widetilde{\Theta},\widetilde{U}\big) dt$$\normalsize is a
solution of system~(\ref{eq:NewtonOperator}). This expression can be
computed at precision~$j$ by a procedure called
\ProcName{ConstantsVariation} in figure~\ref{fig:Algorithm}. \medskip

\subsection{Algorithm}
\label{sec:algo}
We summarize our algorithm in figure~\ref{fig:Algorithm}. This is a
simplified presentation where the technical details are neglected. \medskip
\begin{figure}[!t]
  \small
\begin{center} 
\caption{\label{fig:Algorithm} Local Algebraic Observability Test}
\fbox{\begin{minipage}{\textwidth} \begin{tabbing}
      \Cindent\=\Cindent\=\Cindent\=\Cindent\=\Cindent\=\Cindent\=\Cindent\=\kill
      \CInput :~${\dot{X} - F(X,\Theta,U)}$,~${\ Y - G(X,\Theta,U)}$ \\
      \COutput:~Succeed, a boolean \\[\medskipamount]
      \textbf{Preprocessing} \>\>\>\>\> Construction of the \SLP\ 
      coding $\diff{P}{\dot{X}}, \diff{P}{X}, \diff{P}{\Theta}, \nabla
      P, \Phi_{\Theta}$.
      \\[\medskipamount]
      \textbf{Initialization}
      \>\>\>\>\> Choice of a prime number; \\[\medskipamount]
      \>\>\>\>\> $U$ \Caf Random Power Series mod $t^{n+\ell+1}$;\\
      \>\>\>\>\> Succeed\Caf true; Order\Caf $1$; $\Theta$\Caf Random
      Integers;
      \\
      \>\>\>\>\> $\Lambda$ \Caf $0_{n\times \ell}$; $\Gamma$ \Caf
      $\mathrm{Id}_{n\times n}$; $X$ \Caf Random Integers;
      \\[\medskipamount]
      \Cwhile Order ${\leq n+\ell+1}$ \Cdo \\ \> $W$ \Caf
      $\textup{\ProcName{HomogeneousResolution}} \left(
        \diff{P}{\dot{X}}(\Phi,\Theta)\,\dot{W} +
        \diff{P}{X}(\dot{X},\Phi,\Theta)\,W = 0 \right)\!\! \bmod
      t^{\Order}$;
      \\[\medskipamount]
      \> $(\Phi,\Lambda,\Gamma)$ \Caf $(\Phi,\Lambda,\Gamma) +
      \textup{\ProcName{ConstantsVariation}} \Bigl(W,\nabla P
      \left(\Phi,\Gamma,\Lambda\right)\Bigr) \!\bmod t^{\Order}$;
      \\[\medskipamount]
      \> Increase Order; (Order \Caf $2$ Order);\\
      \Cend \Cwhile
      \\[\medskipamount]
      \> JacobianMatrix \Caf $\textup{\ProcName{Coeffs}}\!\left(
        \nabla Y \left(\Phi,\Gamma,\Lambda\right),
        t^{j},j=0,\dots,n+\ell\right)$; \\[\medskipamount] \CTest
      \>\>\>\>\> \Cif ${n+\ell > \textup{\ProcName{Rank}(JacobianMatrix)}}$ \\
      \>\>\>\>\>\>\Cthen Succeed := false\\
      \>\>\>\>\> \Cend \Cif
\end{tabbing}\end{minipage}}
\end{center}
\normalsize \vspace{-\bigskipamount}
\end{figure}

A preprocessing is necessary to construct, from a \SLP\ 
coding~$\System$, another \SLP\ which encodes the associated linear
variational system~$\nabla P$ and the expressions used during its
integration. This step relies mainly on Theorem~\ref{prop:SLPProp}.
\medskip

The next part of the algorithm consists in the computation at
order~${n+\ell+1}$ of the power series solution of~$\nabla P$.  We
recall that in one iteration, the number of correct coefficients is
doubled (see Theorem~2 in~\cite{Geddes:1979}). \medskip

After the main loop, the procedure \ProcName{Coeffs} evaluates the
\SLP~$\nabla Y$ on the series~$\Phi_{j}$,~$\Gamma_{j}$
and~$\Lambda_{j}$ where~${j = \ln_{2}(n+\ell+1)}$; this furnishes the
coefficients of the jacobian matrix (see
Section~\ref{sec:JacobianMatrixComputation}). \medskip

Last, the rank computations described in
Corollary~\ref{cor:RankCondition} are performed to solve the local
observability problem. \smallskip

If there is more than one output variable, the evaluation of~$\nabla
Y$ and the rank computations which are necessary to determine~$\phi$
can be done in the main loop: the computation can be stopped when the
expected rank is reached or when the computed ranks become stationary.
Thus, we can determine the order of derivation~$\nu$ and avoid useless
computations. \medskip

We now present a rough upper bound for the arithmetic complexity.

\subsection{Arithmetic Complexity Estimation}
\label{sec:Complexity}

\noindent\textbf{Notations:} Hereafter, let~$L$ denote the complexity
of evaluation of the system~$\System$ and let $M(j)$ represent the
multiplication complexity of two series at order~${j+1}$.  Using
classical multiplication formula, we have~${M(j) \in
  \BigO\big(j^{2}\big)}$.

Furthermore, let~$\MatMul{j}$ denotes the number of arithmetic
operations sufficient for the multiplication of two square~${j\times j}$
matrices.  Using classical algorithms, we
have~${\MatMul{n}\in\BigO\big(j^{3}\big)}$.
\begin{proposition}
  The number of arithmetic operations on the base field used in the
  algorithm presented in Section~\ref{sec:algo} is bounded
  by$$\BigO\Biggl(\! M(\nu)\Bigl(\MatMul{n+\ell} + (n+m)L\Bigr) +
  (n+\ell+1)\MatMul{n+\ell}\frac{m\nu}{n+\ell} \!\Biggr)$$
\end{proposition}
\textbf{Proof:} From construction done in
Section~\ref{sec:JacobianMatrixComputation} and
Theorem~\ref{prop:SLPProp}, we conclude that the complexity of
evaluation of the \SLP\ coding~${\ldiff{P}{(\dot{X}, X, \Theta)},
  \nabla P}$ and~${\nabla Y}$ is bounded
by~${\BigO\big(\MatMul{n+\ell} + (n+m)L\big)}$.  Hence, at each step,
the number of arithmetic operations necessary to evaluate this \SLP\ 
on power series truncated at order~$j$, is bounded
by~${\BigO\big(M(j)(\MatMul{n+\ell} + (n+m)L)\big)}$.

Furthermore, the determination of the first~$j$ terms of the solution
series of a system of linear ODE~(\ref{eq:NewtonOperator}) requires
${\BigO\big(M(j)(\MatMul{n} + \MatMul{n+\ell})\big)}$ arithmetic
operations by the well-known method of integrating factors (see {\S}~5.2
in~\cite{BrentKung:1978} for more details).  So, as~${M(j) +
  M\big(\lfloor j/2\rfloor\big) + \dots = \BigO\big(M(j)\big)}$ and as
our Newton operator is quadratic, the arithmetic complexity of the
computations of the jacobian matrix~${\ldiff{(\mathcal{L}^{i}G)_{0
      \leq i \leq \nu} } {(X,\Theta)}}$ is bounded
by~${\BigO\big(M(\nu)\big(\MatMul{n+\ell} + (n + m)L\big)\big)}$.

To conclude, we notice that the cost of a rank computation for
a~${i\times j}$ matrix is $\BigO\big(\MatMul{i}j/i\big)$ if~${i\leq j}$ (see
page~108 in~\cite{BiniPan:1994}).  The
Corollary~\ref{cor:RankCondition} describes the rank computations done
at the end of the main loop of our algorithm.\endproof

\noindent\textbf{Remark:} The specialization of input variables on a
randomly chosen polynomial of degree~${n+\ell}$ increases the evaluation
complexity~$L$ of the system~$\System$ but it does not change the
general complexity of the algorithm.  When the system is not
observable we assume that the index~$\nu$ is~${n+\ell}$ (see
Definition~4 in~\cite{CampbellGear:1995}).  \medskip

We have presented the complexity of our algorithm in term of
arithmetic operations on~$\mathbb{Q}$.  Such an operation requires a
time, roughly, proportional to the size of its operands. Using modular
techniques, we control the growth of the integers involved in the
computations. We estimate now an upper bound on these integers; this
bound will be used in Section~\ref{sec:ProbabilisticAspect} in order
to estimate the probability of success of our algorithm.

\subsection{Growth of the Integers}
\label{subsubsec:Bound}

The forthcoming estimations relies on the formal definition of the
jacobian matrix~${\ldiff{(\mathcal{L}^{i}G)_{0 \leq i \leq \nu}}
  {(X,\Theta)}}$ and are not dependent of the computations described
in Section~\ref{sec:JacobianMatrixComputation} and~\ref{sec:Newton}.

Let us introduce a measure for the size of a~${(n+\ell+r)}$-variate
polynomial which influence the growth of the integers
(see~\cite{CastroHageleMoraisPardo:1999} for more details).
\begin{definition}
  Let~$\mathcal{A}$ be a finite set of non zero integers.  The
  (logarithmic) \emph{height} of~$\mathcal{A}$ is defined
  as~${\height{\mathcal{A}} := \ln |\mathcal{A}|}$
  with~${|\mathcal{A}|:=\max\{|\alpha|+1, \alpha \in \mathcal{A}\}}$.
  
  The height of a polynomial with integer coefficients is defined by
  the height of its set of coefficients.
\end{definition}

We summarize in the following lemma some basic properties of height:
\begin{lemma}
\label{lem:height}
Let~$p_{1}, \dots,p_{s}$ be~${(n+\ell+r)}$-variate polynomials with
integer coefficients,~$x$ an integer and~$\partial$ a partial
derivation ($\ldiff{}{x}$ for example).
\begin{itemize} 
\item ${\height{\partial p} \leq \height{p} + \ln \deg p}$;
\item $\height{p(x)} \leq \height{x} \deg p + \height{p}$;
\item $ht\big(\sum^{s}_{i=1}p_{i}\big) \leq
  \max_{i=1..s}\height{p_{i}} + \ln s$; \item $\height{p_{1}p_{2}}
  \leq \height{p_{1}} + \height{p_{2}} + \min\{\deg p_{1},\deg p_{2}\}
  \ln (n+\ell+r+1)$.
\end{itemize}
\end{lemma}

We use the notations introduced in Section~\ref{sec:Notations} and we
denote by~$h$ (resp.~$d$) the maximum height (resp. degree) of the
numerator and of the denominator of the expression involving in
system~$\System$.
\begin{proposition}
  \label{prop:IntegersGrowth}
  Let~$h_{0}$ be the maximum of heights of the
  integers~$X_{0}$,~$\widetilde{\Theta}$ and of the integer
  coefficients of~$\widetilde{U}$.
  \begin{itemize}
  \item $\height{\denom \mathcal{L}^{j}G\left(X_{0}\right)} \leq \!
    (2j+1) (n+m) \Bigl(\!\bigl(2 \ln (n+\ell+r+1) + h_{0}\bigr)d + h
    \Bigr)$;
  \item $\height{\numer \mathcal{L}^{j}G(X_{0})} \leq \!\!
        \begin{array}{c}
          (2j+1) (n+m) \bigl( (2 \ln (n+\ell+r+1) + h_{0})d + h
          \bigr)  \\
         +\: (j+1) \ln 2n(n+m)d + (2j+1) \ln (2j+1). 
        \end{array}
        $
  \end{itemize}
\end{proposition}
\textbf{Proof:} As we are interested in an upper bound, we do not
consider the reduced form of the fractions~$f_{i}$ and~$g_{i}$
involved in~$\mathcal{L}g$ but we consider that all these fractions
share the same denominator~$q$.  So,~${\mathcal{L} = \big(\sum
  f_{i}\partial_{i}\big)/q}$ and~$q$ is the common denominator of
all~$g_{i}$.

Thus, the degree of these non-reduced numerators and denominators is
bounded by~${(n+m)d}$ and the height by~${(n+m)\big(h + d \ln
  (n+\ell+r+1)\big)}$.  Let us notice that the denominator
of~$\mathcal{L}^{j}g$ is~$q^{2j+1}$; these facts and \mbox{Lemma
  \ref{lem:height}} prove the first part of our proposition.

We prove the second part by induction; let us consider~$(v_{j})_{j \in
  \mathbb{N}}$ the sequence of polynomials defined by the numerator
\mbox{of $g$} as initial condition~$v_{0}$ and by the recurrence
relation~${v_{j+1}:= \sum f_{i}\big(q\partial_{i} v_{j} -
  (2j+1)v_{j}\partial_{i}q\big)}$.  By construction,~$v_{j}$ is equal
to the numerator of~$\mathcal{L}^{j}g$. Thus, the degree of~$v_{j}$ is
bounded \mbox{by $(2j+1)(n+m)d-j$} and we obtain the following
recurrence relation from Lemma~\ref{lem:height}:
$$
\height{v_{j+1}} \leq 2(n+m) \big(2 d \ln (n+\ell+r+1)+h \big) +
\height{v_{j}} + \ln 2n(2j+1)(n+m)d.
$$
This is sufficient to conclude.\endproof

\noindent\textbf{Remark:} 
the use of non-reduced fractions simplifies the previous proof but it
increases the upper bound by a factor~${(n+m)}$ which is not
significant in this presentation.

We have showed that the size of the coefficients of the final
specialized jacobian matrix is mainly linear in the differentiation
index~$\nu$.  But some intermediate computations can require integers
of bigger size.  In order to construct a practical and efficient
algorithm, we have to avoid this growth using modular techniques.
\medskip

Almost all the operations used in our algorithm commute with the
canonical homomorphism from~$\mathbb{Q}$ to a finite
field~$\mathbb{F}_{p}$.  But, when we choose a prime number~$p$, we
have to avoid the cancellation of~$\ldiff{P}{\dot{X}}$ mod~$t$ and of
the determinant of~$\ldiff{(\mathcal{L}^{i}G)_{0 \leq i \leq \nu}}
{(X,\Theta)}$.

The cancellation of~$\ldiff{P}{\dot{X}}$ mod~$t$ can be checked at the
begining of our algorithm.  Thus, the probabilistic aspects concern
mainly the choice of specialization and of a prime number such that
the determinant of~$\ldiff{(\mathcal{L}^{i}G)_{0 \leq i \leq \nu}}
{(X,\Theta)}$ does not vanish modulo~$p$ if this matrix is of full
generic rank.

\subsection{Probabilistic Aspects}
\label{sec:ProbabilisticAspect}

Hereafter, we call a \emph{bad point}, a set of
specializations~${\{X_{0}, \widetilde{\Theta}, \widetilde{U}\}}$ where
the jacobian matrix~${\ldiff{(\mathcal{L}^{i}G)_{0\leq i\leq\nu}}
  {(X,\Theta)}}$ is not of full generic rank. Thus, a bad point is a
zero of the polynomial associated with a minor of this matrix.  We
estimate the probability for a specializations to be a bad point with
the following proposition.
\begin{proposition}[\covek{R. Zippel}~\& \covek{J. Schwartz}~\cite{Zippel:1993}]
\label{prop:ZippelSchwartz} Let~$q$ be a $s$-variate polynomial of
total degree~$D$ and~$\Omega$ a set of integers.  The worst case bound
for the probability that a point in~$\Omega^{s}$ will be a zero of~$q$
is~$D/\card{\Omega}$.
\end{proposition}
This result shows the relation between the choice of the size~$h_{0}$
of the used specializations and the probability of success of our
algorithm.  In fact, as the determinant
of~${\ldiff{(\mathcal{L}^{i}G)_{0 \leq i \leq \nu}} {(X,\Theta)}}$ is
a polynomial of degree bounded by~${D:=(n+\ell)(2\nu+1)(n+m)d}$, a
point in~${\{0,\dots,\mu _{1}D\}^{(n+\ell)(r+1)}}$ is not a bad point
with probability at least~${1-1/\mu _{1}}$.  \medskip

Furthermore, we can estimate the probability that the determinant is
divisible by a prime number~$p$ with the following arithmetic analogue
of Proposition~\ref{prop:ZippelSchwartz}.
\begin{proposition}[{\S}~18 in~\cite{GathenGerhard:1999}]
  For any integers~$a$ and~$b$ such that~${b<a<c}$, the probability
  that a prime number~$p$ between~$b+1$ and~$2b$ divides~$a$ is
  bounded by~$2\ln{c}/b$.
\end{proposition}
From Proposition~\ref{prop:IntegersGrowth} and Lemma~\ref{lem:height},
we can estimate the size of the coefficients of the specialization of
the jacobian matrix~${\ldiff{(\mathcal{L}^{i} G)_{0 \leq i \leq \nu}}
  {(X,\Theta)}}$.  Thus, using Hadamard's inequality, we find the
following rough upper bound for the size of the specialized
determinant:$$
\height{c} := \big(2 \ln (n+\ell+r+1) + h_{0}\big)D +
(n+\ell)(2 \nu +1)\big( (n+m)h + \ln 2nD \big)
$$
Thus, if the computations are performed modulo a prime number~$p$
greater or equal to~$2\height{c}\mu _{2}$ then the probability that the
specialized determinant is not divisible by~$p$ is at
least~${1-1/\mu _{2}}$.  These results lead to the following estimation.
\begin{proposition}
  Let~$\mu $ be an arbitrary positive integer and
  \small$$\begin{array}{ccl}
    D &:=&(n+\ell)(2\nu+1)(n+m)d,\\
    \height{c} &:=& \big(2 \ln(n+\ell+r+1) + \ln D\big)D +
    (n+\ell)(2\nu+1)\big( (n+m)h + \ln 2nD \big).
  \end{array}$$\normalsize
  
  If the matrix~${\ldiff{(\mathcal{L}^{i}G)_{0\leq i\leq\nu}}
    {(X,\Theta)}}$ is of full generic rank then the determinant of
  this matrix specialized on random integers in~${\{0,\dots,\mu D\}}$ is
  not divisible by a prime number~${p > 2\height{c}\mu }$ with
  probability at least~${(1-1/\mu )^{2}}$.
\end{proposition}

\section{Experimental Results}
\label{sec:Conclusion}

We present now some benchmarks from an implementation in
Maple~\cite{SedoglavicURL} of our algorithm.  The Maple computer
algebra system provides almost all the necessary tools to handle the
canonical isomorphism between polynomials and polynomial functions:
this explains why we have chosen it to implement our algorithm. \medskip

The computations summarized in figure~\ref{fig:bench} have been
performed on a personal computer Pentium~III (633~Mhz) with~128Mb of
memory running Linux~2.2 and Maple~V.5. This computer was provided by
the UMS MEDICIS~\cite{MEDICISURL}. \smallskip

\begin{figure}[!t]
  \begin{center}
   \caption{\label{fig:bench} Some benchmarks}    
   {\sf
\begin{tabular}{cccccccl}
  \hline
 System  &$m$&$\nu$&$\ell$&$n$ &$r$& $L$ & time in s.\\ \hline
 V1987   & 2 &  8  &   5  &  4 &   & 17 & $\quad \ \: 0.8$\\  
 R1986   & 2 & 14  &   9  &  4 & 1 & 19 & $\quad \ \: 1.5$\\  
 MV1991  & 2 & 14  &   8  &  5 & 2 & 59 & $\quad \ \: 2.4$\\  
 MW2000  & 3 & 18  &  14  &  4 &   & 67 & $\quad \ \: 5.7$\\  
 KD1999  & 2 & 19  &  14  &  5 & 2 & 34 & $\quad \ \: 6.$\\    
 G1995   & 1 & 23  &  17  &  5 &   & 46 & $\quad 10.$\\        
 SHH1997 & 1 & 23  &  13  &  9 &   & 38 & $\quad 13.5$\\\hline
\end{tabular}
}
  \end{center} \vspace{-\bigskipamount}
\end{figure}

These results show that the index of differentiation is a significant
characteristic of the complexity of algorithm presented in
Section~\ref{sec:algo}.  Furthermore, the last example of the array
shows that the complexity of evaluation have a significative influence
and that the total number of multiplications is clearly less
significant than the number of multiplications between state and input
variables.

\subsection{Certifying the result}
\label{sec:group}

As shown in Corollary~\ref{cor:RankCondition}, the local observability
property is associated to the fact that the jacobian matrix is of full
rank.  Our algorithm computes the generic rank of this matrix. When it
is maximal, the result is certainly correct.  Hence, if this algorithm
states that a model is observable then this result is certified (it is
a RP-complexity class test, see {\S}~25.8 in~\cite{GathenGerhard:1999}).
\medskip

If there is a non empty set~${O\subset X \cup\Theta}$ of non observable
variables and parameters, the observable parameters can be randomly
specialized and there is an infinitesimal transformation acting on the
non observable state variables and parameters,
$$
\mathcal{S} := \sum_{x\in O\cap X} s_{x} \diff{\;\;}{x} + \sum_{\theta
  \in O\cap\Theta} s_{\theta} \diff{\;\;}{\theta},
$$
which leaves invariant the outputs~$G$ and the vector field
associated to the model. This leads to the following linear system of
PDE's:
$$
\left\lbrace
\begin{array}{cl}
\big[\mathcal{S},\mathcal{L} \big] & =0,\\
\mathcal{S} G &= 0.
\end{array}\right.
$$
This system of PDE can be difficult to solve; nevertheless, we are
not interested in the whole Lie algebra but in any non trivial
subalgebra which can certified our result.

Furthermore, our algorithm decreases the number of unknown of the
original problem.  Hence, in many cases of practical interest, there
is a rather straightforward solution (compare
with~\cite{Ollivier:1999}). For example, these computations have been
performed in less than a hour with Maple for the following examples.

\subsection{Examples}
\label{app:sys}

We present now the examples indicated in figure~\ref{fig:bench}, the
answer of our algorithm and some results of the method sketched in
section~\ref{sec:group}.  We just give the non observable parameters
and variables; the other one are observable.

\subsubsection*{V1987 \ Model of a flow reactor to pyrolyze methane}
This example is taken from~\cite{Vajda:1987}.
$$
\left\lbrace
  \begin{array}{ccl}
    \dot{x}_{1} &=& -x_{1}(k_{1} + k_{2}x_{4}) + k_{5}x_{3}x_{4}, \\
    \dot{x}_{2} &=& k_{2}x_{1}x_{4} - (k_{3}+k_{4})x_{2},\\
    \dot{x}_{3} &=& k_{4}x_{2} - k_{5}x_{3}x_{4},\\
    \dot{x}_{4} &=& x_{1}(k_{1} + k_{2}x_{4}) + 2k_{3}x_{2} - k_{5}x_{3}x_{4},
    \\[\smallskipamount]
    y_{1}  &=& x_{1},\\
    y_{2}  &=& x_{2}.
  \end{array}
\right.
$$
Our Maple implementation certifies that all the variables
and the parameters are observable.

\subsubsection*{R1986 \ A pharmacokinetic model}
This example is taken from~\cite{Raksanyi:1986}. The letter~$u$
denotes an input.
$$
\left\lbrace
  \begin{array}{ccl}
    \dot{x}_{1} &=& u - (c_{1} + c_{2})x_{1}, \\
    \dot{x}_{2} &=& c_{1}x_{1}-(c_{3} + c_{6} + c_{7})x_{2} + c_{5}x_{4}, \\
    \dot{x}_{3} &=& c_{2}x_{1} + c_{3}x_{2} - c_{4}x_{3}, \\
    \dot{x}_{4} &=& c_{6}x_{2} - c_{5}x_{4}, \\[\smallskipamount]
    y_{1} &=& c_{8}x_{3}, \\
    y_{2} &=& c_{9}x_{2}.
  \end{array}
\right.
$$
Our Maple implementation gives the following results:
\begin{itemize}
\item the variables~${\{x_{2},x_{3},x_{4}\}}$ and the
parameters~${\{c_{1},c_{2},c_{3},c_{7},c_{8},c_{9}\}}$ are not
observable;
\item the transcendence degree of the field extension~${k\langle U,Y
    \rangle \hookrightarrow \mathcal{K}}$ is~$1$.
\end{itemize} 
Further computations show that the following one parameter group
$$
\begin{array}{ll}
\begin{array}{ccc}
x_{2} & \rightarrow & \lambda x_{2} \\[\smallskipamount]
x_{3} & \rightarrow & ((1-\lambda)c_{1} + c_{2}) 
                       {x_{3}}/{c_{2}} \\[\smallskipamount]
x_{4} & \rightarrow & \lambda x_{4}  \\[\smallskipamount] 
c_{1} & \rightarrow & \lambda c_{1} \\[\smallskipamount]
c_{2} & \rightarrow & (1 - \lambda)c_{1} + c_{2} 
\end{array} &
\begin{array}{ccc}
c_{3} & \rightarrow & \big((1-\lambda)c_{1} + c_{2}\big){c_{3}}/{\lambda c_{2}}
                      \\[\medskipamount]
c_{7} & \rightarrow & c_{7} - {c_{3}(c_{1} + c_{2})(1-\lambda)}/
                      {\lambda c_{2}} \\[\medskipamount]
c_{8} & \rightarrow & -{c_{8}c_{2}}/{\big((1-\lambda)c_{1} + c_{2}\big)} 
                      \\[\medskipamount]
c_{9} & \rightarrow & {c_{9}}/{\lambda} 
\end{array}
\end{array}
$$
is composed of symmetries which leave the vector field and the
output invariant.

\subsubsection*{MV1991 \ Model for an induction motor}
This example is taken from~\cite{MarinoValigi:1991}.  The
letters~$u_{x}$ and~$u_{y}$ denote inputs.
$$
\left\lbrace
  \begin{array}{cclccl}
    \sigma &=& L_{s} - \frac{M^{2}}{L_{r}}, &
    \gamma_{N} & = & \frac{M^{2} R_{r} + {L_{r}}^{2}R_{s}}{\sigma {L_{r}}^{2}},
    \\[\bigskipamount]

     \dot{\omega} &=& \multicolumn{4}{l}{\frac{n_{p} M}{J L_{r}}(\Psi_{x}I_{y}
                        - \Psi_{y}I_{x})- \frac{T_{L}}{J},} \\[\medskipamount]

     \dot{\Psi}_{x} &=& \multicolumn{4}{l}{- \frac{R_{r}}{L_{r}} \Psi_{x} 
                        - n_{p} \omega \Psi_{y} + \frac{R_{r}}{L_{r}} M I_{x},}
                         \\[\medskipamount]

     \dot{\Psi}_{y} &=& \multicolumn{4}{l}{n_{p}\omega\Psi_{x}-\frac{R_{r}}
                        {L_{r}} \Psi_{y} + \frac{R_{r}}{L_{r}} M I_{y},} 
                        \\[\medskipamount]

     \dot{I}_{x} & = & \multicolumn{4}{l}{\frac{M R_{r}}{\sigma {L_{r}}^{2}}
                        \Psi_{x} + \frac{n_{p} M}{\sigma L_{r}}\omega\Psi_{y} 
                        - \gamma_{N} I_{x} + \frac{u_{x}}{\sigma},}
                        \\[\medskipamount]

     \dot{I}_{y} & = & \multicolumn{4}{l}{- \frac{n_{p} M}{\sigma L_{r}}\omega
                        \Psi_{x} + \frac{M R_{r}}{\sigma {L_{r}}^{2}}\Psi_{y} 
                        - \gamma_{N} I_{y} + \frac{u_{y}}{\sigma},}
                        \\[\bigskipamount]

      y_{1} &=& \omega, \\
      y_{2} &=& {\Psi_{x}}^{2} + {\Psi_{y}}^{2}.                       
  \end{array}
\right.
$$ 
Our Maple implementation gives the following results:
\begin{itemize}
\item the variables~${\{I_{x}, I_{y}\}} $ and the
  parameters~${\{M,L_{s},R_{s},L_{r},R_{r},J,T_{l}\}}$ are not
  observable;
\item the transcendence degree of the field extension~${k\langle U,Y
    \rangle \hookrightarrow \mathcal{K}}$ is~$1$.
\end{itemize} 
Further computations show that the following one parameter group
$$
{\{I_{x}, I_{y}, M, L_{s}, R_{s}, L_{r}, R_{r}, J, T_{l}\}}
\rightarrow
{\{\lambda I_{x}, \lambda I_{y}, M\lambda, L_{s}\lambda, R_{s}\lambda, 
L_{r}/\lambda, \lambda R_{r}, \lambda J, \lambda T_{l}\}}
$$
is composed of symmetries which leave the vector field and the
output invariant.

\subsubsection*{MW2000 Multispecies model for the transmission of pathogens}
This example is taken from~\cite{MargariaWhite:2000}.
$$
\left\lbrace
  \begin{array}{ccl}
b & = & \mu + c_{1}(y_{1}+y_{12}), \\          
\lambda_{1} &=& \beta_{1}(y_{1}+y_{12}),         \\
\lambda_{2} &=& \beta_{2}(y_{2}+y_{12}) + I_{2},  \\[\medskipamount]
\dot{x}_{12} &=& (1-\theta_{1}-\theta_{2})b - (m_{1}\lambda_{1} +
m_{2}\lambda_{2} + \mu )x_{12} +\\
&&  (\nu_{1}+\tau)y_{1} + (\nu_{2}+\tau)y_{2} + \tau y_{12}, \\

\dot{y}_{1} &=& \theta_{1}b + m_{1}\lambda_{1}x_{12} + \nu_{2}y_{12}
-\big((1-\pi_{2})m_{2}\lambda_{2}+\nu_{1}+\mu +c_{1}+\tau\big)y_{1}, \\

\dot{y}_{2} &=& \theta_{2}b + m_{2}\lambda_{2}x_{12} + \nu_{1} y_{12}
-\big((1-\pi_{1})m_{1}\lambda_{1}+\nu_{2}+\mu +\tau\big)y_{2}, \\

\dot{y}_{12} &=& (1-\pi_{1})m_{1}\lambda_{1}y_{2} +
(1-\pi_{2})m_{2}\lambda_{2}y_{1} - (\nu_{1}+\nu_{2}+\mu +c_{1}+\tau)y_{12}, 
\\[\medskipamount]

o_{1} &=& x_{12}+ y_{1} + y_{2} + y_{12}, \\
o_{2} &=& y_{1} + y_{12}, \\
o_{3} &=& y_{2} + y_{12}.
  \end{array}
\right.
$$
Our Maple implementation gives the following results:
\begin{itemize}
\item with the exception of~${\{\beta_{1}, \beta_{2}, I_{2}, m_{1},
    m_{2}\}}$, all the parameters are observable;
\item the transcendence degree of the field extension~${k\langle U,Y
    \rangle \hookrightarrow \mathcal{K}}$ is~$2$.
\end{itemize} 
Further computations show that the following two parameters group
$$
{\{\beta_{1},\ \beta_{2},\ I_{2},\ m_{1},\ m_{2}\} \rightarrow
  \{\beta_{1}/l_{1},\ \beta_{2}/l_{2},\ I_{2}/l_{2},\ l_{1} m_{1},\ l_{2}
  m_{2}\}}
$$
is composed of symmetries which leave the vector field and the
output invariant. \smallskip

Let us notice that the output~$o_{1}$ is in fact a constraint equal
to~$1$. Hence, our model can be composed of relations of order zero
which can be considered as supplementary outputs.

\subsubsection*{KD1999 \ Model for a chemical reactor} 
This example is taken from~\cite{KumarDaoutidis:1999}.
$$
\left\lbrace
      \begin{array}{ccl} 
        \dot{C}_{A} & = &\frac{F_{A}}{V}(C_{A0}-C_{A}) 
                        - k_{0}e^{-E/RT}C_{A},\\[\smallskipamount]
        \dot{C}_{B} & = & -\frac{F_{A}}{V}C_{B}) 
                        + k_{0}e^{-E/RT}C_{A},\\[\smallskipamount]

        \dot{T} & = & \frac{F_{A}}{V}(T_{A}-T) - k_{0}e^{-E/RT}C_{A}
                        \frac{\Delta H_{r}}{\rho c_{p}} + \frac{U}{\rho c_{p}}
                        \frac{T_{j}-T}{V},\\[\smallskipamount]
        \dot{T}_{j} & = & \frac{F_{h}}{V_{h}}(T_{h}-T_{j}) -
                         \frac{U}{\rho_{h} c_{ph}}
                        \frac{T_{j}-T}{V_{h}},\\[\medskipamount]
        y_{1} & = & C_{B}, \\
        y_{2} & = & T.
\end{array} \right.
$$
We denote by A the Arrhenius' law~$e^{-E/RT}$ and we add the
ordinary differential equation~${\dot{A} = E A \dot{T}/(R T^{2})}$ to
the model.  Our Maple implementation gives the following results:
\begin{itemize}
\item the variable~$A$ and the parameters~${\{E,R,\Delta H_{r},U,
\rho, c_{p}, \rho_{h}, c_{ph},k_{0}\}}$ are not observable;
\item the transcendence degree of the field extension~${k\langle
    U,Y \rangle \hookrightarrow \mathcal{K}}$ is~$5$.
\end{itemize} 
Further computations show that the following five parameters group
$$
\begin{array}{ccc}
\begin{array}{ccc}
A  & \rightarrow & \lambda_{1}A \\
k_{0} & \rightarrow & k_{0}/\lambda_{1} \\
E & \rightarrow & \lambda_{2}E \\
R  & \rightarrow & \lambda_{2}R
\end{array} &
\begin{array}{ccc}
\rho & \rightarrow & \lambda_{3}\rho\\
c_{p} & \rightarrow & \lambda_{4}c_{p} \\[\smallskipamount]
\Delta H_{r}  & \rightarrow & \lambda_{3}\lambda_{4}\Delta H_{r} 
\end{array} &
\begin{array}{ccc}
U & \rightarrow & \lambda_{3}\lambda_{4}U\\
c_{ph}  & \rightarrow & \lambda_{5}c_{ph}  \\
\rho_{h} & \rightarrow & \lambda_{3}\lambda_{4}\rho_{h}/\lambda_{5}
\end{array}
\end{array}
$$
is composed of symmetries which leave the vector field and the
output invariant.

\subsubsection*{G1995 \ Model of Circadian oscillations in the Drosophila 
  period protein} This example is described in introduction.
\subsubsection*{SHH1997 \ Model of a part of the blood coagulation mechanism}

This example is taken from~\cite{StortelderHemkerHemker:1997}.
$$
\begin{array}{lll}
r_{1} = \frac{kc_{X}X \cdot RVV}{km_{X}+X}, &
r_{2} = ki_{Xa}Xa, &
r_{3} = \frac{kc_{V}V \cdot IIa}{km_{V}+V}, \\[\smallskipamount]
r_{4} = k_{PT}Va \cdot Xa \cdot PL, &
r_{5} = k_{PL} PT, &
r_{6} = \frac{kc_{II} \cdot II \cdot PT}{km_{II}+II}, \\[\smallskipamount]
r_{7} = \frac{kc_{2} II \cdot Xa}{km_{2}+II}, &
r_{8} = ki_{IIa\alpha_{2}M} \cdot IIa, &
r_{9} = ki_{IIaATIII} \cdot IIa.
\end{array}
$$
$$\left\lbrace \begin{array}{ccl}
    \dot{X} &=&  -r_{1}, \\
    \dot{Xa}&=& r_{1} - r_{2} - r_{4} + r_{5},        \\
    \dot{V} &=&-r_{3},        \\
    \dot{Va}&=&r_{3} - r_{4} + r_{5},         \\
    \dot{PL} &=&- r_{4} + r_{5},           \\
    \dot{PT} &=&r_{4} - r_{5},        \\
    \dot{II} &=&- r_{6} - r_{7},        \\
    \dot{IIa}   &=& r_{6} + r_{7} - r_{8} - r_{9},       \\
    \dot{IIa\alpha_{2}M} &=& r_{9},       \\[\smallskipamount]
    y &=& IIa + \frac{556}{1000}IIa\alpha_{2}M.
\end{array} \right.$$ 
Our Maple implementation gives the following results:
\begin{itemize} 
\item the parameters~${\{kc_{X}, km_{X}, kc_{V}, km_{V}, k_{PT},
    kc_{II}, kc_{2}\}}$ and \\ the variables~${\{X, Xa, V, Va, PL,
    PT\}}$ are not observable; 
\item the transcendence degree of the field
  extension~${k\langle U,Y \rangle \hookrightarrow \mathcal{K}}$
  is~$1$.
\end{itemize} 
Further computations show that the following one parameter group
$$
\begin{array}{ccc}
\begin{array}{ccc}
X  & \rightarrow & \lambda X \\
Xa & \rightarrow & \lambda Xa \\
V  & \rightarrow & \lambda V \\
Va & \rightarrow & \lambda Va
\end{array} &
\begin{array}{ccc}
PL & \rightarrow & \lambda PL\\
PT & \rightarrow & \lambda PT \\
kc_{X}  & \rightarrow & \lambda kc_{X} \\
km_{X} & \rightarrow & \lambda km_{X}
\end{array} &
\begin{array}{ccc}
kc_{V}  & \rightarrow & \lambda kc_{V}  \\
km_{V} & \rightarrow & \lambda km_{V} \\
k_{PT}  & \rightarrow & k_{PT}/\lambda^{2}   \\
kc_{II} & \rightarrow & kc_{II}/\lambda  \\
kc_{2} & \rightarrow & kc_{2}/\lambda 
\end{array}
\end{array}
$$
is composed of symmetries which leave the vector field and the
output invariant.
\paragraph{Acknowledgment:} It is a pleasure to thank \covek{M. Giusti},
\covek{G. Lecerf}, \covek{F. Ollivier} and \covek{{\'E}. Schost} for their
contributions, helpful comments and illuminating discussions.

\bibliographystyle{acm} \bibliography{$HOME/Biblio/references}
\end{document}